\documentclass[12pt,a4paper,twoside]{article}
\usepackage{umlaut}
\usepackage{a4wide}
\usepackage[T1]{fontenc}    
\usepackage{ae,aecompl}     
\usepackage{amsmath}        
\usepackage{amssymb}
\usepackage{amsthm}
\usepackage{amsfonts}
\usepackage{amscd}
\usepackage{bbm}            
\usepackage{bm}             
\usepackage{url}            
\usepackage{paralist}       
\usepackage{xspace}      
\newtheorem{thm}{Theorem}[section]

\newtheorem{prop}[thm]{Proposition}
\theoremstyle{definition}

\newtheorem{defn}[thm]{Definition}
\newtheorem{expl}[thm]{Example}

\usepackage{knesermacros}

\begin{document}

\title{On generalized {K}neser hypergraph colorings}
\author{Carsten E. M. C. Lange\thanks{Supported by a Postdoc Stipend of the
    Mittag-Leffler Institute}\\[-3pt]
\small University of Washington\\[-5pt]
\small Department of Mathematics\\[-5pt]
\small Box 354350\\[-5pt]
\small Seattle, WA 98195-4350\\[-5pt]
\small \url{lange@math.tu-berlin.de}\setcounter{footnote}{6}%
\and Günter M. Ziegler\thanks{Partially supported by Deutsche
  Forschungs-Gemeinschaft (DFG), via Project~B9 of the
  \textsc{Matheon} Research Center ``Mathematics for Key
  Technologies'' (FZT86), the Research Group ``Algorithms, Structure,
  Randomness'' (Project ZI~475/3), and a Leibniz grant}\\[-3pt]
\small TU Berlin\\[-5pt]
\small Inst.\ Mathematics, MA 6-2\\[-5pt]
\small D-10623 Berlin, Germany\\[-5pt]
\small \url{ziegler@math.tu-berlin.de}} \date{\today}
\date{February 2006}
\maketitle

\begin{abstract}
\noindent
In Ziegler (2002), the second author presented a lower bound for the
chromatic numbers of hypergraphs $\KG{r}{\pmb s}{\calS}$,
``generalized $r$-uniform Kneser hypergraphs with intersection
multiplicities $\pmb s$.'' It generalized previous lower bounds by 
K\v{r}\'{\i}\v{z} (1992/2000) for the case ${\pmb s}=(1,\dots,1)$
without intersection multiplicities, and by Sarkaria (1990) for
$\calS=\tbinom{[n]}k$.
Here we discuss subtleties and difficulties that arise
for intersection multiplicities $s_i>1$:
\begin{compactenum}[1.]
\item In the presence of intersection multiplicities, there are two
  different versions of a ``Kneser hypergraph,'' depending on whether
  one admits hypergraph edges that are multisets rather than sets. We
  show that the chromatic numbers are substantially different for the
  two concepts of hypergraphs.  The lower bounds of Sarkaria (1990)
  and Ziegler (2002) apply only to the multiset version.
\item The reductions to the case of prime~$r$ in the proofs by Sarkaria and
  by Ziegler work only if the intersection multiplicities are strictly smaller
  than the largest prime factor of~$r$. Currently we have no valid
  proof for the lower bound result in the other cases.
\end{compactenum}
We also show that all uniform hypergraphs 
without multiset edges can be represented 
as generalized Kneser hypergraphs.
\end{abstract}

\section{Introduction} \label{sec_intro}

The ``generalized Kneser hypergraphs with intersection
multiplicities,'' as studied in \cite{Z77}, arose from the graphs
(implicitly) studied by Kneser \cite{Kneser} in several subsequent
generalization steps. They do form a large class of hypergraphs ---
indeed, we will show in Section~\ref{sec_universality} that every
uniform hypergraph without multiplicities can be represented in
this model.

When writing \cite{Z77}, the second author overlooked that the 
edges of the generalized Kneser hypergraph $\KG r{\pmb s}\calS$ 
with intersection multiplicities $\pmb s=(s_1,\dots,s_n)$
could be multisets if $s_i>1$: Their edges can have repeated
elements, as pointed out in~\cite{lange:phd_thesis}. Thus 
$\KG r{\pmb s}\calS$ is not a hypergraph in the
traditional sense of Berge~\cite{Berge}, where the edges have to be
sets.  If one does not allow for repeated elements in the edges, then
this yields a sub-hypergraph without multiplicities,
\[
\kg r{\pmb s}(\calS)\ \ \subseteq\ \ \KG r{\pmb s}(\calS).
\]
Both for $\kg r{\pmb s}(\calS)$ and for $\KG r{\pmb s}(\calS)$ we are
faced with the problem to determine the chromatic number: How many
colors are needed for the sets in $\calS$ if monochromatic hypergraph
edges are forbidden?  Clearly we have
$\chi(\kg r{\pmb s}\calS)\le\chi( \KG r{\pmb s}\calS)$,
but the two values can be far apart, as we will see below.
\medskip

\noindent
In this note, we discuss the chromatic numbers of generalized
Kneser hypergraphs with intersection multiplicities,
in view of some main topics and results from~\cite{Z77}. This
includes errata and clarifications announced in~\cite{Z77erratum}:
\smallskip

\begin{compactdesc}
\item[{\cite[Lemma~3.1]{Z77}}] described an explicit coloring for the
  special case of $\calS=\tbinom{[n]}k$ and constant 
  ${\pmb s}=(s,...,s)$. 
  This coloring is valid for
  generalized Kneser hypergraphs \emph{with} multiplicies, which
  yields
\[
\chi(\kg rs\tbinom nk)\ \ \le\ \ \chi( \KG rs\tbinom nk)\ \ \le\ \ 
1+\big\lceil\tfrac1{\lfloor\tfrac{r-1}s\rfloor}\tfrac{ns-rk+1}s\big\rceil.
\]

\item[{\cite[Theorem~5.1]{Z77}}] states a lower bound 
   \[
\lceil \tfrac{1}{r-1} \ \cd{r}{\pmb s}{\calS}\rceil
\ \ \le\ \ \chr{\KG{r}{\pmb s}{\calS}}
   \]
for the generalized Kneser hypergraphs. This lower bound holds for 
generalized Kneser hypergraphs \emph{with} multiplicities and is
not valid for $\chr{\kg{r}{\pmb s}{\calS}}$ as we 
will see in Section~\ref{sec_examples}.

Moreover, Karsten Vogel (Magdeburg) has noticed that the reduction
to the prime case in the proof of \cite[Theorem~5.1]{Z77} fails
when the intersection multiplicities are not
smaller than the prime factors of~$r$.
As we will analyze in Section~\ref{sec_induction},
we get only the following theorem (with a combinatorial proof):

\begin{thm}\label{thm:main_thm}
  Let $\calS \subseteq 2^{[n]}$ be a set family, and let
  the intersection multiplicities $s_i\ge1$ be smaller than the 
  largest prime factor of~$r\ge2$.  Then
   \[
     \chr{\KG{r}{\pmb s}{\calS}} \ge \lceil \tfrac{1}{r-1} \ \cd{r}{\pmb s}{\calS}\rceil.
   \]
\end{thm}

In particular, the conditions of this theorem are satisfied 
\begin{compactitem}[$\bullet$]
\item if there are no intersection multiplicities, 
  $\pmb s=(1,...,1)$. In this case 
  $\kg{r}{(1,...,1)}\calS=\KG{r}{(1,...,1)}\calS$, and
  Theorem~\ref{thm:main_thm} reduces to the main result of
  K\v{r}\'{\i}\v{z} \cite{Kri1,Kri2}, and
  \item in the case when $r$ is prime, for arbitrary $s_i<r$.
\end{compactitem}
We still believe that the theorem is valid for arbitrary $s_i<r$,
as stated by \cite[Thm.~5.1]{Z77},
but we have no proof for this generality --- not even for the 
``complete uniform''
case when $\calS=\tbinom{[n]}k$; indeed, we will see in
Section~\ref{sec_induction} that the induction proof by
Sarkaria~\cite{Sarkaria-kneser} for this case  is not valid, either.

\item[{\cite[Example~7.2]{Z77}}] analyzed some Kneser hypergraphs
  \emph{without} multiplicities, including a computation of
  $\chi(\kg4{(2,...,2)}\tbinom{[n]}2)$.
  Thus in Section~\ref{sec_examples} we discuss
  families of hypergraphs that include $\KG4{(2,...,2)}\tbinom{[n]}2$.
  We also collect evidence towards the conjecture that the upper bound
$\chi( \KG r{(s,...,s)}\tbinom nk)\le
1+\big\lceil\tfrac1{\lfloor\tfrac{r-1}s\rfloor}\tfrac{ns-rk+1}s\big\rceil$
is tight in general.
\end{compactdesc}

\section{Preliminaries} \label{sec_prelims}

In this section we review the fundamental concepts for this study;
compare \cite[Sect.~2]{Z77}.  Let $n\ge1$ and denote
$[n]:=\{1,\dots,n\}$.  By $\pmb s$ we denote a vector of positive
integers $\pmb s=(s_1,\dots,s_n)$.  Throughout $r\ge2$ denotes an integer.  
We write~$\pmb s < r$ if $s_{i} < r$ for all~$i$.

\begin{defn}[$\pmb s$-disjoint sets]
Subsets $S_1,\dots,S_r$ of~$[n]$ are $\pmb s$\emph{-disjoint} if each
$i \in [n]$ occurs in at most~$s_i$ of the sets~$S_k$. 
Note that here equalities $S_k= S_\ell$ are possible.
\end{defn}

To illustrate this definition consider $n=3$ and $\pmb s = (3,2,1)$.
The subsets $\{1,2\}$, $\{1,2\}$, and $\{2,3\}$ are not $\pmb s$-disjoint 
because $2$ occurs in all three sets. On the other hand, $\{1,2\}$, $\{1,2\}$, 
and $\{1,3\}$ are $\pmb s$-disjoint, and so are $\{1,2\}$, $\{1,3\}$, and $\{2\}$.

\begin{defn}[$\pmb s$-disjoint $\pmb r$-colorability defect 
  {\cite[p.~673]{Z77}}] Let $[n]^{\pmb s}$ denote the multiset in
  which the element~$i\in[n]$ occurs with multiplicity~$s_{i}$.  We
  denote the cardinality of $[n]^{\pmb s}$ counting multiplicities
  by~$\overline n$.
  
  The $\pmb s$\emph{-disjoint} $r$\emph{-colorability defect}
  $\cd{r}{\pmb s}{\calS}$ of a set $\calS \subseteq 2^{[n]}$ is the
  minimal number of elements one has to remove from the
  multiset~$[n]^{\pmb s}$ such that the remaining multiset can be
  covered by~$r$ subsets of~$[n]$ that do not contain any
  element from~$\calS$:
\[
  \cd{r}{\pmb s}{\calS} = {\overline n} -
                             \max \set{\sum_{j=1}^{r}|R_{j}|}
{\begin{matrix}
  R_{1}, {\ldots} , R_{r} \subseteq [n]\ \pmb s\text{-disjoint subsets}\\
  \text{and } S \not \subseteq R_{j} 
  \text{ for all } S \in \calS \text{ and all }j
 \end{matrix}}.
\]
\end{defn}

Here the sets $R_\ell$ need not be distinct, and they may be empty.
Note for further reference that $\cd{r}{\pmb s}{\varnothing}>0$
if $s_{i} > r$ for some~$i$.

As example we consider $n=3$, $\pmb s=(3,2,1)$, and
${\cal S} =\{ S_{1} \}$ where $S_{1}=\{2,3\}$. We are therefore never allowed
to pick $\{1,2,3\}$ or $\{2,3\}$ as one of the sets $R_{j}$.
For $r=1$, $R_{1}=\{1,2\}$ is a possible choice of largest 
cardinality, thus $\cd{1}{(3,2,1)}{\calS} = 6-2=4$. For $r=2$, we may 
pick $R_{1}=\{1,2\}$ and $R_{2}=\{1,3\}$ as examples that maximize~$|R_{1}|+|R_{2}|$, 
hence $\cd{2}{(3,2,1)}{\calS} =2$. For $r=3$, the value $|R_1|+|R_2|+|R_3|$
is maximized by  $R_{1}=\{1,2\}$, $R_{2}=\{1,2\}$, and $R_{3}=\{1,3\}$, so
$\cd{3}{(3,2,1)}{\calS} =0$.

\begin{defn}[$\pmb r$-uniform hypergraphs with/without multiplicities]
  An $r$-\emph{multi\-subset} $X$ of~$[n]$ is an unordered collection 
  of $r$ elements $x_1,\dots,x_r$ of $[n]$ that need not be distinct.
  We denote it by $X=\{\!\{ x_{1},{\ldots} ,x_{r} \}\!\}$.
  
  An \emph{$r$-uniform hypergraph} in the sense of Berge \cite{Berge}
  (without multiplicities) is a pair~$(V,E)$ that consists of a
  finite set of \emph{vertices} $V$ and a set $E$ of \emph{edges},
  which are $r$-subsets of~$V$.
  
  An \emph{$r$-uniform hypergraph with multiplicities} is a
  pair~$(V,E)$ that consists of a finite set of \emph{vertices} $V$
  and a set $E$ of \emph{edges}, which are $r$-multisubsets of~$V$.

  A hypergraph (with multiplicities) is \emph{loop-free} if every
  edge has at least two distinct elements (vertices).
In the following, all hypergraphs are supposed to be loop-free. 
\end{defn}

According to this definition, hypergraphs do not have multiple
edges; this makes sense for our purposes, since multiple
edges are irrelevant for coloring.

The loop-free $r$-uniform hypergraph analog of the complete graph
$K_{n}$ on $n$ vertices is denoted by $K^{r}_{n}$. The vertex set of
$K^{r}_{n}$ is $[n]$, and the edges are all the $r$-multisubsets of
$[n]$ that contain at least two distinct vertices. 
Thus $K^r_n$ has $\big(\!\!\tbinom nr\!\!\big)-n=\tbinom{n+r-1}r-n$ edges.
The analogous complete $r$-uniform hypergraph $k^r_n$ without
multiplicities has $\tbinom nr$ edges.

\begin{defn}[$\pmb r$-uniform $\pmb s$-disjoint Kneser hypergraphs]
For any finite set $\calS = \{S_1,\dots ,S_m\}$ of non-empty subsets
of~$[n]$, the $r$-\emph{uniform} $\pmb s$-\emph{disjoint}
\emph{Kneser hypergraph} $\KG{r}{\pmb s}{\calS}$ \emph{with multiplicities} 
has the vertex set $\calS$ and the edge set
\[
  \edges {\KG{r}{\pmb s}{\calS}}\ \ :=\ \ 
\{\text{ all $\pmb s$-disjoint $r$-multisets whose elements are
  sets~$S_i\in \calS$ }\}.
\]
If $s_i<r$ for all $i\in[n]$ then $\KG{r}{\pmb s}{\calS}$ is loop-free.

The $r$\emph{-uniform} $\pmb s$\emph{-disjoint Kneser hypergraph}
$\kg{r}{\pmb s}{\calS}$ \emph{without multiplicities} has the same
vertex set $\calS$, but all of its edges are sets rather than multisets:
\[
    \edges {\kg{r}{\pmb s}{\calS}}\ \ :=\ \ 
\{\text{ all $\pmb s$-disjoint $r$-subsets of~$\calS$ }\}.
\]
\end{defn}

The generalized Kneser hypergraphs $\kg rs\calS$ are loop-free
for any $\pmb s$.

We use $\KG rs\calS$ as a shorthand for $\KG r{(s,...,s)}\calS$
in the case of constant intersection multiplicity
$\pmb s=(s,...,s)$, and 
similarly we write $\kg rs\calS$ and $\cd rs\calS$.

The previously defined complete $r$-uniform hypergraphs~$K^r_n$ and $k^r_n$ are
examples of $r$-uniform $\pmb s$-disjoint Kneser hypergraphs. We have 
$K^r_n = \KG{r}{r-1}{\tbinom{[n]}{1}}$, $k^r_n = \kg{r}{r-1}{\tbinom{[n]}{1}}$, and
in this particular situation 
$\KG{r}{1}{\tbinom{[n]}{1}}=\kg{r}{r-1}{\tbinom{[n]}{1}}$.

We can obtain
$\kg{r}{\pmb s}{\calS}$ from $\KG{r}{\pmb s}{\calS}$ by discarding
edges. In this sense, $\kg{r}{\pmb s}{\calS}$ is a subhypergraph of
$\KG{r}{\pmb s}{\calS}$.  In the special case that 
$s_i\equiv1$ we have
$\KG{r}{s}{\calS} = \kg{r}{s}{\calS}$ since 
pairwise disjoint non-empty sets are distinct. In particular, 
for $r=2$ and $s_{i} \equiv 1$ both definitions specialize to the
generalized Kneser graph of $\calS \subseteq 2^{[n]}$.
 
\begin{defn}[hypergraph colorings 
 \cite{erdos_hajnal_on_chromatic_number_of_graph_and_set_systems}]
  A coloring of an $r$-uniform hypergraph~$H$
  (multiplicity-free or not) with~$m$ colors is a map $c:\vertices{H}
  \rightarrow [m]$ that assigns to each vertex of~$H$ a color such
  that no edge is monochromatic, that is, for each $e \in \edges{H}$
  we have $| \set{c(x)}{x\in e} | \ge 2$. Any coloring $c$ of
  $H$ by $m$ colors induces a homomorphism $H\rightarrow K^r_m$ of
  hypergraphs. The \emph{chromatic number}~$\chr{H}$ is the smallest
  number~$m$ such that there is a coloring of~$H$ with~$m$ colors.
\end{defn}

\section{How general are generalized Kneser hypergraphs?}
\label{sec_universality}%

Matou\v{s}ek \& Ziegler \cite[p.~76]{Z83} observed that every (finite,
simple) graph can be represented as a Kneser graph: For any $G=(V,E)$
there is a set system $\calS=\set{S_v}{v\in V}\subseteq2^{[m]}$,
for some~$m$,
such that $S_v\cap S_w=\varnothing$ if and only if $\{v,w\}\in E$.
Thus it is natural to ask which hypergraphs (with or without
multiplicities) can be represented as generalized Kneser hypergraphs.
The following proposition collects our answers to this question.

For this, call a hypergraph \emph{up-monotone} if for
$e,e'\in\big(\!\!\tbinom{[n]}r\!\!\big)$ with $e\in E$ we also have
$e'\in E$ whenever the support of~$e'$ contains that of $e$. Every
$r$-uniform hypergraph without multiplicities is up-monotone, as is
every generalized Kneser hypergraph $\KG r{r-1}\calS$.
  
A hypergraph $H=([n],E)$ is \emph{convex} if every
integral weight vector $(a_1,\dots,a_n)$ in the convex hull of
multiplicity vectors of edges of $H$ (thus $0\le a_i<r$) is the
multiplicity vector of an edge of~$H$.
   
\begin{prop}\label{prop:represent_by_Kneser}~
\begin{compactenum}[\rm(1)]
\item There are $r$-uniform hypergraphs without multiplicities that
  cannot be represented as a Kneser hypergraph $\KG r1\calS$.
  
\item An $r$-uniform hypergraph $H=([n],E)$ with multiplicities can be
  represented as $\KG r{r-1}\calS$ if and only if it is up-monotone.
  
  In particular, every $r$-uniform hypergraph without multiplicities
  can be represented as a Kneser hypergraph $\KG r{r-1}\calS$.
  
\item If an $r$-uniform hypergraph is representable by a generalized
  Kneser hypergraph with intersection multiplicities then it is
  convex. (The converse is not true.)
  
  In particular, there are $r$-uniform hypergraphs with multiplicities
  that cannot be represented as a Kneser hypergraph 
  $\KG r{\pmb s}\calS$.
\end{compactenum}
\end{prop}

\begin{proof}
  (1). Consider $([4],\{124,134,234\})$.  If 
  $\KG31\{S_1,\dots,S_4\}$ has $\{S_1,S_2,S_4\}$,
  $\{S_1,S_3,S_4\}$ and $\{S_2,S_3,S_4\}$ as edges, then each of the
  triples of sets is pairwise disjoint, so in particular $S_1,S_2,S_3$
  are pairwise disjoint. Thus also $\{S_1,S_2,S_3\}$ is an edge in the
  Kneser hypergraph.%
\smallskip
  
(2). The following construction generalizes the construction
  for graphs in~\cite{Z83}. Let
  $H=([n],E)$ be up-monotone, and let $\bar H=([n],\bar E)$ be the
  complementary hypergraph of $H$, i.e.\ the hypergraph that has the
  same vertices as~$H$ and all edges of $K^{r}_{n}$ that are not edges
  of~$H$.
  Define the set system $\calS=\set{S_i}{i\in[n]}$ by
\[ S_i\ \ :=\ \ \{i\}\cup\set{\bar e\in\bar E}{i\in\bar e}. \]
  The $S_i$ are clearly distinct.
  If $e=\{\!\{ i_1,\dots,i_r\}\!\}$ is an edge of~$H$, then 
\[ S_{i_1}\cap\dots\cap S_{i_r}\ \ =\ \ 
   \{i_1\}\cap\dots\cap\{i_r\}\ \cap\ \set{\bar e\in\bar
     E}{i_1,\dots,i_r\in\bar e},
\]
  where the first part is empty since $H$ does not have loops (so the
  $i_k$ cannot all be equal) and the last set is empty since $H$ is
  up-monotone.  Thus $S_{i_1}\cap\dots\cap S_{i_r}$ is an edge of the
  Kneser hypergraph. Conversely, if 
  $S_{i_1}\cap\dots\cap S_{i_r}=\varnothing$, then in particular 
  it does not contain the element $e=\{i_1\dots,i_r\}$, so $e\in E$.  
\smallskip
  
(3). The intersection multiplicities $s_i$
  define the hypergraph $\KG r{\pmb s}\calS$ as 
  a subgraph of~$K^r_m$, for $m=|\calS|$, by \emph{linear} conditions 
  on the multiplicity vectors of the edges. 
  
  For an example consider $([3],\{113,223\})$.  
  If $\KG3{\pmb s}\{S_1,S_2,S_3\}$, with $S_1,S_2,S_3\subseteq[m]$, 
    does not have $\{S_1,S_2,S_3\}$ as an edge, then there is some
  $i\in[m]$ such that $S_1,S_2,S_3$ contain~$i$ more than $s_i$~times.
  However, that cannot be if both $\{\!\{S_1,S_1,S_3\}\!\}$ and
  $\{\!\{S_2,S_2,S_3\}\!\}$ are edges, so $S_1,S_1,S_3$ and
  $S_2,S_2,S_3$ contain $i$ at most~$s_i$ times.

  An example of a convex uniform hypergraph with multiplicities
  that cannot be represented as a generalized Kneser hypergraph is 
  $([3],\{112,223\})$.
\end{proof}

For the purpose of coloring, any hypergraph $H=(V,E)$ with
multiplicities can be replaced by an up-monotone uniform hypergraph
with multiplicities, on the same ground set, and with the same
chromatic number: For this replace each edge $e$ by all multisets of
cardinality $r$ which contain the support of~$e$, for
some large enough~$r$.
By Proposition~\ref{prop:represent_by_Kneser}(2),
the resulting $r$-uniform hypergraph with multiplicities $H'$ can be
represented as a generalized Kneser graph, which 
yields topological lower bounds for $\chr H$ in terms
of the colorability defect of~$H'$. 
In particular, this applies to (non-uniform) hypergraphs
in the sense of Berge.

\section{Two counterexamples} \label{sec_examples}

The purpose of this section will be to show that the lower bound 
$\chi(\KG r{\pmb s}\calS)\ge\lceil\frac1{r-1}\cd rs\calS\rceil$
of Theorem~\ref{thm:main_thm} is not valid for $\kg{r}{s}{\calS}$.
Indeed, the proof for \cite[Thm.~5.1]{Z77} is valid only
if multiplicities are included: The argument at \cite[p.~679]{Z77}
yields $p$ subsets $S_1,\dots,S_p$ of~$[n]$ that 
are $\pmb s$-disjoint, 
but they need not be pairwise different. 

\begin{expl}
For $n\ge5$ and $r\ge4$ with $n\ge r-1$, let
$\calS:=\{ 12, 13, \dots , 1n,\, 23, 45 \}\subset\tbinom{[n]}2$.
Then all edges of~$\kg{r}{r-2}{\calS}$ are of the form
$\{1i_1,\dots,1i_{r-2},23,45\}$, so they contain both
$23$ and $45$. 
Thus $\chr {\kg{r}{r-2}{\calS}} = 2$.
A straightforward argument 
(cf.~\cite[p.~83]{lange:phd_thesis} for details) 
shows that $\cd{r}{r-2}{\calS} = 3r - 10$. For $r > 8$ this yields
$\cd{r}{r-2}{\calS} > (r-1)\chr{\kg{r}{r-2}{\calS}}$.
\end{expl}

The next example shows that for Kneser hypergraphs without
multiplicities, the colorability defect lower bound does not even hold
in the special case of $\calS=\tbinom{[n]}k$.  
(Sarkaria~\cite{Sarkaria-kneser} speaks of ``$p$-tuples of
$S$-subsets'', so his treatment clearly concerns the
Kneser hypergraphs $\KG rs\tbinom{[n]}k$
\emph{with} multisets as hypergraph edges.)

\begin{expl}
For $n,r\ge4$ the hypergraph $\kg r{r-1}\tbinom{[n]}2$ has a
greedy $(n-2)$-coloring, by $c:S\mapsto\min\{\min S,n-2\}$.
(The hypergraph is non-empty if $r\le \tbinom n2$. 
Its chromatic number will be computed in Example~\ref{example:kg}.)

On the other hand, 
$\cd{r}{r-1}{\tbinom{[n]}{2}} = \max\{n(r-1) - r,0\}=(n-1)(r-1)-1$
by~\cite[Lemma~3.2]{Z77}. Thus
\[
  (r-1)\chr {\kg  r{r-1}{\tbinom{[n]}2}}\ \ <\ \ 
            \cd r{r-1}{\tbinom{[n]}2}.
\]
\end{expl}

\section{The induction to non-prime cases} \label{sec_induction}

For the case of Theorem~\ref{thm:main_thm} when $p$ is a prime
Ziegler~\cite{Z77} has given a combinatorial proof;
an alternative topological proof was given by
Lange~\cite[Sect.~4.4]{lange:phd_thesis}. 
The special case when $\calS=\tbinom{[n]}k$
is due to Sarkaria~\cite{Sarkaria-kneser}.

In this section, we show that the reductions of the situation with
general $r$ to the case of prime~$r$ 
by Sarkaria \cite{Sarkaria-kneser} and by Ziegler \cite{Z77} are both
incomplete. We also argue that argument given in \cite{Z77} suffices
to establish the result in the generality of
Theorem~\ref{thm:main_thm}.  \medskip

Sarkaria's proof \cite[(3.2)]{Sarkaria-kneser} starts with the
assumption that $\KG p{j-1}\tbinom{[N]}S$ has an $M$-coloring, with
$(p-1)M < N(j-1) - p(S-1) = \cd{p}{j-1}{\tbinom{[N]}{S}}$,
for some non-prime $p=p_{1}p_{2}$. Then one constructs a 
coloring of $\KG {p_{1}}{j-1}{\tbinom{[N]}{N'}}$ with 
\[
  N' := M(p_{2}-1) + p_{2}(S-1) + 1,
\]
and tries to get a monochromatic $(j-1)$-disjoint $p_1$-tuple of
$N'$-subsets of~$[N]$. The argument fails if $N'$ is larger than~$N$,
so there won't be \emph{any} $N'$-subsets of~$[N]$.
Concrete parameters where this happens are
$p=4$,
$p_1=p_2=2$,
$j=4$,
$S=2$,
$M=N-2$, 
which yields 
$N'=M+3=N+1$.
(The problem does not occur for $j\le3$, so in particular 
the proof specializes correctly to the case $j=2$ treated by
Alon, Frankl \& Lov\'asz~\cite{AlonFranklLovasz}.)
\medskip

Ziegler's reduction to the prime case in \cite[pp.~679-680]{Z77}
is an extension of K{\v r}{\' i}{\v z}' proof~\cite{Kri2}, which in
turn generalizes the argument of Alon, Frankl and Lov\'asz. 
Let $r=r'r''$ with $r' \le r''$. The goal is to 
derive a contradiction if we assume that 
$\cd{r}{\pmb s}{\calS} > (r-1)\chr{\KG{r}{\pmb s}{\calS}}$. A
crucial ingredient is the set
\[
  \calT\ \ :=\ \ \set{N \subseteq [n]}
  {\cd{r'}{1}{\calS_{|N}} > (r'-1)\chr{\KG{r}{\pmb s}{\calS}}}
  \label{eqn:defn_of_U}
\]
where $\calS_{|N}$ denotes the elements of $\calS$ that are subsets
of~$N$. One then wants to argue that
\[
(r''-1)\chr{\KG{r''}{\pmb s}{\calT}}
\ \ \ge\ \         \cd{r''}{\pmb s}{\calT}
\]
But this can be concluded by induction only if $\pmb s < r''$.
Moreover, it definitely fails if $s_{i} > r''$ for some $i$ and
$\calT=\varnothing$: 
In this situation $\chr{\KG{r''}{\pmb s}{\varnothing}} = 0$ since
there are no vertices to color, but 
$\cd{r''}{\pmb s}{\varnothing}>0$ since at least $s_{i}-r''$ elements have to be
removed removed from $[n]^{\pmb s}$ to cover the remaining elements
with $r''$ subsets of~$[n]$).
The case $\calT=\varnothing$ can occur, as we have seen above
for the special case of $\calS=\tbinom{[N]}S$.

Thus, \cite[Thm.~5.1]{Z77} can currently only be established
in the generality given above as Theorem~\ref{thm:main_thm}.
To establish
this, one uses the induction given at \cite[pp.~679-680]{Z77},
factoring non-prime $r=r'r''$ so that $r''$ is the largest prime
number that divides~$r$.

\section{More Examples} 
\label{sec_special_case}

In \cite[Sect.~7]{Z77} the second author had raised the question 
whether the upper bound of~\cite[Lemma~3.1]{Z77}
\begin{equation}\tag{$*$}
  \label{eq:upper_bound}
  \textstyle
   \chr{\KG rs}{\tbinom{[n]}k}  \ \ \le\ \ 
   1 + \Big\lceil \frac{1}{\left\lfloor \frac{r-1}{s} \right\rfloor}\frac{ns-rk+1}{s}\Big\rceil
\end{equation}
is always tight, for $n\ge k\ge2$, $r>s\ge2$, $rk\le sn$.
In \cite[Example~7.2]{Z77} he had claimed that
(\ref{eq:upper_bound}) is not sharp for $\KG42{\tbinom{[n]}2}$.
However, this is not true: The analysis given there referred to 
the corresponding Kneser hypergraph without
multiplicities, that is, it established that
\[\textstyle
\chi(\kg42{\tbinom{[n]}2})\ \ =\ \ 
n-\big\lfloor\sqrt{2n+\tfrac14}-\tfrac12\big\rfloor.
\]
Thus the tightness question is open for now.  By
Theorem~\ref{thm:main_thm}, (\ref{eq:upper_bound}) is tight if $s$ is
smaller than the largest prime factor of~$r$, and divides~$r-1$
(cf.~\cite[Cor.~7.1]{Z77}).  The following example yields more cases
where (\ref{eq:upper_bound}) is tight, including the case
of~$\KG42{\tbinom{[n]}2}$. 

\begin{expl}\label{example:KG}
Assume that $k=2$ and $\lfloor \tfrac{r-1}{s} \rfloor = 1$,
i.e.\ $\tfrac{r}2 \le s < r-1$. Then 
\[\textstyle
    \chr{\KG rs{\tbinom{[n]}2}} \ \ =\ \ 
   1 + \big\lceil \frac{ns-2r+1}{s}\big\rceil
\ =\    1 + n - \big\lfloor \frac{2r-1}{s}\big\rfloor.
\]
Indeed, the vertices of $H=\KG{r}{s}{\tbinom{[n]}{2}}$ are the edges
of a complete graph~$K_n$. 
By $s\ge\lceil\tfrac r2\rceil$, an edge of~$H$
cannot contain two disjoint edges from~$K_n$.
Thus the edges of~$H$ are supported only on stars 
or on triangles --- the latter is permitted if $s\ge\tfrac{2r}3$.
Thus the possible color classes $C\subset E(K_n)$ are
stars, or they are triangles --- the latter is permitted if
$s<\tfrac{2r}3$. In either
case the greedy colorings that provide the upper bound
are optimal: $n-1$ colors are needed for $\tfrac{2r}{3}\le s \le r-1$,
while $n-2$ colors are needed for $\tfrac r2\le s<\tfrac{2r}3$.
\end{expl}

\begin{expl}\label{example:kg}
The Kneser hypergraphs without multiplicities
$\kg r{r-1}{\tbinom{[n]}2}$ have chromatic numbers
\[\chr{\kg r{r-1}{\tbinom{[n]}2}}\ \ =\ \ 
\begin{cases}
 \lceil \tfrac{1}{r-1}\tbinom{n}2 \rceil & n < r, \\
  n-\lfloor\tfrac r2\rfloor              & 2 \le r \le n.
\end{cases}
\]
Indeed, any edge of this hypergraph forms an $r$-set of edges in~$K_n$
that is not a star. Thus for a color class we can use any star, or any
set of at most $r-1$ edges.  An optimal coloring in case of $2 \le r \le n$ 
uses $n-r$ stars, and then 
$\lceil\tfrac1{r-1}\tbinom r2\rceil=\lceil\tfrac r2\rceil$ edge
sets of size at most $r-1$ to cover the remaining uncolored subgraph
$K_r$. If $n < r$, an optimal colouring uses 
$\lceil\tfrac1{r-1}\tbinom n2\rceil$ sets of size at most $r-1$.
\end{expl}

In summary, we see that
\begin{eqnarray*}
  n-\lfloor\tfrac r2\rfloor\ \ =\ \ 
\chr{\kg r{r-1}{\tbinom{[n]}2}} &\ll&
\chr{\KG r{r-1}{\tbinom{[n]}2}} \ \ =\ \ n-1 \\[-7pt]
\noalign{\noindent and}
n-\big\lfloor\sqrt{2n+\tfrac14}-\tfrac12\big\rfloor\ \ =\ \ 
\chr{\kg 42{\tbinom{[n]}2}} &\ll&
\chr{\KG 42{\tbinom{[n]}2}} \ \ =\ \ n-1 
\end{eqnarray*}
for sufficiently large $n$ and $r$. This shows a
huge difference between the chromatic numbers of 
generalized Kneser hypergraphs with and without multiplicities.
\begin{small}

\end{small}

\end{document}